\documentclass[11pt]{amsart}
\usepackage{amscd}
\usepackage{amsfonts}
\usepackage{amssymb,latexsym}
\usepackage{enumerate}
\usepackage{amsmath}
\usepackage{mathrsfs}
\usepackage{amsthm}
\usepackage{hyperref}
\usepackage[square, comma, sort&compress, numbers]{natbib}
\usepackage{bm}
\usepackage{esint}
\usepackage{fancyhdr}
\usepackage{lastpage}
\usepackage{layout}
\usepackage{ulem}
\usepackage{extarrows}
\usepackage{geometry}
\usepackage{comment}
\newcommand{\ol}{\overline}

\newcommand{\pa}{\partial}

\DeclareMathOperator\tr{tr}
\DeclareMathOperator\Id{id}
\DeclareMathOperator\rank{rank}
\DeclareMathOperator\Ric{Ric}

\DeclareMathOperator\dist{dist}

\DeclareMathOperator\End{End}

\makeatletter
\@namedef{subjclassname@1991}{Mathematics Subject Classification}
\makeatother
\begin{document}
	\newcounter{remark}
	\newcounter{theor}
	\setcounter{remark}{0}
	\setcounter{theor}{1}
	\newtheorem{claim}{Claim}[section]
	\newtheorem{theorem}{Theorem}[section]
	\newtheorem{lemma}{Lemma}[section]
	\newtheorem{corollary}{Corollary}[section]
	\newtheorem{corollarys}{Corollary}
	\newtheorem{question}{Question}
	\newtheorem{proposition}{Proposition}[section]
	\newtheorem{defn}{Definition}[section]
	\newtheorem{examp}{Example}[section]
	\newtheorem{assumption}{Assumption}[section]
	\newtheorem{rem}{Remark}[section]
	\newtheorem*{theorem1}{Theorem}
	\newtheorem*{maintheorem}{Main Theorem}
	\newtheorem*{maintheorem2}{Main Theorem \uppercase\expandafter{\romannumeral2}}
	\numberwithin{equation}{section}
	\title[Difference estimates for Chern connections on homolomorphic vector bundles]
	{Difference estimates for Chern connections on homolomorphic vector bundles}
	\author{Chenghong He}
	\address{Chenghong He, School of Mathematics and Statistics, Nanjing University of Science and Technology, Nanjing 210094, People's Republic of China}
	\email{Hech@njust.edu.cn}
	\author{Di Wu}
	\address{Di Wu, School of Mathematics and Statistics, Nanjing University of Science and Technology, Nanjing 210094, People's Republic of China}
	\email{wudi@njust.edu.cn}
	\author{Xi Zhang}
	\address{Xi Zhang, School of Mathematics and Statistics, Nanjing University of Science and Technology, Nanjing 210094, People's Republic of China}
	\email{mathzx@njust.edu.cn}
	\thanks{The authors are supported by the National Natural Science Foundation of China No.12601105, No.12431004, No.12141104, the China Postdoctoral Science Foundation 2025T181124, 2023M731699, the Jiangsu Funding Program for Excellent Postdoctoral Talent 2022ZB282, the Natural Science Foundation of Jiangsu Province No.BK20230902, the Fundamental Research Fund for Central Universities WZCT202601009, and the Stability Support Program for Basic Research WZJC202601001.}
	\maketitle
	\begin{abstract}
		This paper investigates the uniqueness issue of Chern connections from the Cheng-Yau type estimate point of view. 
		We establish local difference estimates for Chern connections on holomorphic vector bundles over complete noncompact K\"{a}hler manifolds with Ricci curvature bounded from below. These local estimates not only enable us to obtain uniqueness results for Hermitian-Einstein connections on complete noncompact K\"{a}hler manifolds with nonnegative Ricci curvature, but also apply to study the Liouville theorem for complex Monge-Amp\`{e}re equation on $\mathbb{C}^n$ as well as the Chern flatness and K\"{a}hlerness for Hermitian metrics on complete noncompact Ricci flat K\"{a}hler manifolds.
	\end{abstract}
	
	\section{Introduction}
	Let $(E,\bar\partial_E)$ be a holomorphic vector bundle over a K\"ahler manifold. 
	A Hermitian metric $H$ on the underlying complex vector bundle is called Hermitian-Einstein if 
	\begin{equation}\begin{split}\label{HE}
			\Lambda F_H^\perp\triangleq\Lambda F_H-\frac{\tr\Lambda F_H}{\rank(E)}\Id_E=0,
	\end{split}\end{equation}	
	where $\Lambda$ denotes the contraction with the K\"{a}hler form 
	and $F_H$ is the curvature of Chern connection $D_H$. 
	We call such $D_H$ a Hermitian-Einstein connection and sometimes it is also referred to as a Hermitian-Yang-Mills connection in the literature. 
	Thanks to a series of works including Narasimhan-Seshadri \cite{NS1965}, Mehta-Ramanathan \cite{MR1984}, 
	Donaldson \cite{Donaldson85,Donaldson87} and Uhlenbeck-Yau \cite{UhlenbeckYau86}, 
	the celebrated Donaldson-Uhlenbeck-Yau theorem on compact K\"ahler manifolds indicates 
	the existence of Hermitian-Einstein metrics on Mumford-Takemoto poly-stable bundles. 
	Furthermore, it is not hard to see two Hermitian-Einstein metrics with same determinant on a Mumford-Takemoto stable bundle  must be the same.

	\par
	There have been several results solving the Hermitian-Einstein equation $(\ref{HE})$ on noncompact K\"{a}hler manifolds, see \cite{Ba1993,MochizukiHE,MS2010,Ni2002,NiRen,Si1988,WuZhangHE,ZhangJ} and references therein.  
	However, the uniqueness issue seems more subtle on noncompact manifolds.  
	To our knowledge, there are not many known results in this direction. Even for the most important model that the underlying space can be compactified, the uniqueness issue was unknown until the recent work \cite{MochizukiHE} of Mochizuki, while the existence part was established by Simpson \cite{Si1988} more than thirty years ago.
	One may also consult \cite{WuZhangHE} for an alternative approach.

	\par 
We shall study the uniqueness issue of Hermitian-Einstein problem from a new point of view. 
	A fundamental source of inspiration is the classical Cheng-Yau gradient estimate.
	\begin{theorem}[\cite{ChengYau75,Yau75}]\label{thm:CY}
		Let $(M,g)$ be a complete
		Riemannian manifold with dimension $n$, $\Ric_g\geq -\kappa g$ and $\kappa\geq0$ on a geodesic ball $B_{R}$.
		If $u\in C^\infty(B_R)$ is a positive harmonic function, 
		then there is a positive constant $C$ depending only on $\dim M$ such that
		\begin{equation}
			\sup_{B_{R/2}}\frac{|\nabla u|}{u}\leq\frac{C(1+\sqrt{\kappa}R)}{R}.
		\end{equation}
		In particular, there is no nontrivial glocal positive harmonic function if $\Ric_g\geq0$.
	\end{theorem}

	\par
	Cheng-Yau's approach has been a paradigm and an inspiration for more than fifty years. 
	Many important results for other problems are deeply influenced by Theorem \ref{thm:CY}. 
	
	\par
	We start with the simplest case that two Hermitian-Einstein metrics $H,K$ are mutually bounded, 
	conformal and satisfy $\tr\Lambda F_H=\tr\Lambda F_K$. 
	It follows directly from Theorem \ref{thm:CY}
	that they have the same Chern connection provided the base K\"{a}hler manifold is complete and has nonnegative Ricci curvature. Naively, this naturally motivates us to study the general uniqueness issue from the philosophy of Cheng-Yau estimate. 
	We expect two mutually bounded Hermitian-Einstein metrics with same determinant must have the same Chern connection.
	\par Given two Hermitian metrics $H,K$ on $E$, we define $h=K^{-1}H\in A^0(M,\End(E))$ by
	\begin{equation}\begin{split}\label{relative-endomorphism}
			H(\bullet,\bullet)=K(h(\bullet),\bullet).
	\end{split}\end{equation}
Then $h$ is self-adjoint with respect to $K$ and $H$. Let $\lambda_{R,\max}$ be the supremum of maximum eigenvalue of $h$ on $B_R$ and $\lambda_{R,\min}$ be the infimum of minimal eigenvalue of $h$ on $B_R$.
\par Throughout this paper, unless
otherwise mentioned, inner products and norms are taken with respect to $H$ and $g$. Moreover, we always denote by $C$ a positive constant which only depends on $\dim_{\mathbb{C}}M$ and may vary from line to line.		
	\par 
	Our central idea is to establish a local difference estimate for Chern connections.
		\begin{theorem}\label{thm:main11}
		Let $(E,\bar\partial_E)$ be a holomorphic vector bundle over a complete 
		K\"ahler manifold $(M,g)$ with 
		$\Ric_g\geq-\kappa g$ and $\kappa\geq0$ on a geodesic ball $B_{R}$.
		Suppose $H$, $K$ are two Hermitian metrics on $E$ 
		with $\det H=\det K$. If $H$ is Hermitian-Einstein, then it holds
		\begin{equation}\label{main-estimate-general}
			\begin{aligned}
				\sup_{B_{R/2}}|D_H-D_K|
				&\leq C(A_R\lambda_{R,\max})^{1/2} \frac{1+\sqrt{\kappa}R}{R}\\
				&+C(A_R\lambda_{R,\max})^{1/2}(\sup_{B_R}|F_K^\perp|+
				\sup_{B_R}|\Lambda F_K^\perp|_K)^{1/2}\\
				&+C(A_R\lambda_{R,\max})^{1/3}\sup_{B_R}|D_K(\Lambda F_K^\perp)|^{1/3},
			\end{aligned}
		\end{equation}
		where $A_R=\sup\limits_{B_R}\tr h^{-1}$ and $h=K^{-1}H$.
	\end{theorem}
	\par	
	If $\Ric_g\geq0$ and $K$ is projectively flat (that is $F_K^{\perp}=0$),
	Theorem~\ref{thm:main11} implies
	\begin{equation}\label{estimatereduce}
		\sup_{B_{R/2}}|D_H-D_K|
		\leq\frac{C\sqrt{\rank(E)}}{R}
		(\frac{\lambda_{R,\max}}
		{\lambda_{R,\min}})^{1/2},
	\end{equation}
	where we have used $A_R\lambda_{R,\min}\leq \rank(E)$. 
	Therefore we obtain
	\begin{corollary}\label{cor1}
		Let $(E,\bar\partial_E)$ be a holomorphic vector bundle over a complete noncompact 
		K\"ahler manifold with nonnegative Ricci curvature.
		Suppose $H$, $K$ are two Hermitian metrics on $E$ with $\det H=\det K$. 
		If $H$ is Hermitian-Einstein, $K$ is projectively flat and 
		\begin{equation}
			\frac{\lambda_{R,\max}}
		{\lambda_{R,\min}}=o(R^2),
				\end{equation}
			then $D_H=D_K$.
	\end{corollary}
	In particular, we establish the expected uniqueness of projectively flat Chern connections. 
	One can also interpret Corollary \ref{cor1} as a rigidity result 
	as it yields a sufficient condition for a Hermitian-Einstein connection to be projectively flat. 
    The same rigidity requires the Chern number restriction on
    compact K\"{a}hler manifolds, 
	while Simpson \cite{Si1988} investigated the same issue on the complement of an ample divisor in a compact K\"{a}hler manifold.
	\par 
	Essentially, Cheng-Yau's philosophy is to give an estimate of harmonic objects. It is easy to see that $F_H\in A^{1,1}(M,\End(E))$ is a harmonic section 
	if $H$ is Hermitian-Einstein and $\tr\Lambda F_H$ is constant, 
	therefore a parallel theorem to Theorem \ref{thm:CY} should give an estimate of $F_H$. However, since the difference between two Chern connections is not harmonic and in this sense, Cheng-Yau's philosophy is used in a nonstandard way and it brings difficulty in the local difference estimate with additional terms appearing on the right side. 
	We illustrate in the Appendix that they can not be all dropped by showing there are two mutually bounded Hermitian-Einstein metrics with same determinant 
	but different Chern connections on $\mathbb{C}^2$. 
	\par 
	Although there is a counterexample showing the projectively flat assumption 
	cannot be replaced by the Hermitian-Einstein one, 
	we find a necessary and sufficient condition for the validity of the expected uniqueness
	by the following local difference estimate.
		\begin{theorem}\label{thm:main22}
		Let $(E,\bar\partial_E)$ be a holomorphic vector bundle over a complete 
		K\"ahler manifold $(M,g)$ with 
		$\Ric_g\geq-\kappa g$ and $\kappa\geq0$ on a geodesic ball $B_{R}$.
		Suppose $H$, $K$ are two Hermitian metrics on $E$ 
		with $\det H=\det K$. If $H$ is Hermitian-Einstein, then it holds
		\begin{equation}\label{curvature-action-estimate-general}
			\begin{aligned}
				\sup_{B_{R/2}}|D_H-D_K|
				&\leq C(A_R\lambda_{R,\max})^{1/2}
				\frac{1+\sqrt{\kappa}R}{R}
				\\&+C(A_R\lambda_{R,\max})^{1/2}
				(\sup_{B_R}|h^{-1}[F_K,h]|
				+\sup_{B_R}|\Lambda F_K^\perp|
				+\sup_{B_R}|\Lambda F_K^\perp|_K)^{1/2}
				\\&+C(A_R\lambda_{R,\max})^{1/3}
				(\sup_{B_R}|h^{-1}\Lambda D_K[F_K,h]|
				+\sup_{B_R}|D_K(\Lambda F_K^\perp)|)^{1/3},
			\end{aligned}
		\end{equation}
where $A_R=\sup\limits_{B_R}\tr h^{-1}$ and $h=K^{-1}H$.
	\end{theorem}
	\begin{corollary}\label{cor2}
		Let $(E,\bar\partial_E)$ be a holomorphic vector bundle
		over a complete noncompact K\"ahler manifold with nonnegative Ricci curvature.
		Suppose $H$, $K$ are two Hermitian-Einstein metrics on $E$ with $\det H=\det K$ and	
		\begin{equation}\label{growth-condition-main}
			\frac{\lambda_{R,\max}}
			{\lambda_{R,\min}}=o(R^2).
		\end{equation}
		Then $D_H=D_K$ if and only if $[F_K,h]=0$.
	\end{corollary}
	Below we compare our results with previous works \cite{MochizukiHE,WuZhangHE}. 
	Following \cite{WuZhangHE}, 
	two Hermitian metrics $H,K$ are said to be mutually compatible 
	if 
\begin{equation}\begin{split}
\det H=\det K,\ K^{-1}H\in L^\infty,\ \ol\pa_E(K^{-1}H)\in L^2. 
\end{split}\end{equation}
	Mochizuki \cite{MochizukiHE} established the uniqueness of mutually compatible Hermitian-Einstein metrics 
	on a stable Hermitian holomorphic vector bundle over a Zariski open subset in a compact K\"{a}hler manifold. The proof essentially uses the Uhlenbeck-Yau regularity \cite{UhlenbeckYau86}. 
	Without the recourse to this deep preliminary, 
	Wu-Zhang \cite{WuZhangHE} gave a simpler alternative proof and it applies to more general situations. In \cite{MochizukiHE,WuZhangHE}, the uniqueness of Hermitian-Einstein connections is a direct consequence of the Liouville theorem for subharmonic functions. The key point there is to prove the uniqueness of Hermitian-Einstein metrics 
	from that of Hermitian-Einstein connections under stability condition. 
	Unlike \cite{MochizukiHE,WuZhangHE}, 
	this paper focuses on the uniqueness of Hermitian-Einstein connections themselves. 
\par In the following, we apply our result to
holomorphic tangent bundles .
\par Let $g_0$ be the standard metric on $\mathbb{C}^n$ with K\"{a}hler form $\omega_0$ and trivial connection $\nabla_{g_0}$. We obtain the following Liouville theorem for complex Monge-Amp\`{e}re equation on $\mathbb{C}^n$.
\begin{theorem}\label{thm:main3}
Let $\omega$ be a Ricci flat K\"{a}hler form on $\mathbb{C}^n$ and there are positive functions $c_R,C_R$ such that
\begin{equation}\begin{split}
		c_R\omega_0|_{B_R}\leq\omega|_{B_R}\leq C_R\omega_0|_{B_{R}},\ \frac{C_R}{c_R^2}=o(R^2),\ B_R=\{z\in\mathbb{C}^n,|z|<R\}. 
\end{split}\end{equation}
Then $\nabla_{g_0}\omega=0$ and equivalently, there is a constant Hermitian matrix $A=(A_{ij})$ such that
\begin{equation}\begin{split}
		\omega=\sqrt{-1}\sum\limits_{i,j=1}^nA_{ij}dz^i\wedge d\ol{z}^j.
\end{split}\end{equation}
\end{theorem}
In particular, it recovers the classical Liouville theorem of Riebesehl-Schulz.
\begin{corollary}[\cite{RS1984}]\label{cor:RS}
	Let $\omega$ be a K\"{a}hler form on $\mathbb{C}^n$ and there is a constant $C$ such that
\begin{equation}\begin{split}
		\omega^n=\omega_{0}^n,\
		C^{-1}\omega_0\leq\omega\leq C\omega_0.
\end{split}\end{equation}
Then $\nabla_{g_0}\omega=0$ and equivalently, there is a constant Hermitian matrix $A=(A_{ij})$ such that
\begin{equation}\begin{split}
		\omega=\sqrt{-1}\sum\limits_{i,j=1}^nA_{ij}dz^i\wedge d\ol{z}^j.
\end{split}\end{equation}
\end{corollary}
We denote by $\nabla_{\bullet}$ the Levi-Civita connection of the Riemannian metric $\bullet$. We get
\begin{theorem}\label{thm4}
Let $(M,g)$ be a complete noncompact Ricci flat K\"{a}hler manifold. Suppose $\tilde{g}$ is a Chern flat Hermitian metric and there are positive functions $c_R,C_R$ such that
\begin{equation}\begin{split}
		c_Rg|_{B_R}\leq\tilde{g}|_{B_R}\leq C_Rg|_{B_{R}},\ \frac{C_R}{c_R}=o(R^2),\ B_R=\{q\in M|\dist_g(p,q)<R\},\ p\in M.
\end{split}\end{equation}
Then $\nabla_g=\nabla_{\tilde{g}}$. In particular, $g$ is flat and $\tilde{g}$ is K\"{a}hler.
\end{theorem}
\begin{corollary}\label{cor3}
Let $(M,g)$ be a complete noncompact Ricci flat K\"{a}hler manifold. Suppose $\tilde{g}$ is a Hermitian metric with zero Chern holomorphic sectional curvature and there are positive functions $c_R,C_R$ such that
\begin{equation}\begin{split}
		c_Rg|_{B_R}\leq\tilde{g}|_{B_R}\leq C_Rg|_{B_{R}},\ \frac{C_R}{c_R}=o(R^2),\ B_R=\{q\in M|\dist_g(p,q)<R\},\ p\in M.
\end{split}\end{equation}
Then $\tilde{g}$ is Chern flat if and only if it is K\"{a}hler. If these equivalent conditions hold, $\nabla_g=\nabla_{\tilde{g}}$.
\end{corollary}
Another motivation comes from a conjecture in Hermitian geometry stating that a compact Hermitian
manifold with constant Chern holomorphic sectional curvature $c$ is K\"{a}hler if $c\neq0$ and Chern flat if $c=0$. We are interested in the zero constant case. The conjecture is known to be true in dimension two \cite{ADM1996,Bal1985,BG1985} and dimension three \cite{CZ2026,QT2026}. Thanks to \cite{Bo1958}, compact Chern flat manifolds are quotients of complex Lie groups which may be non-K\"{a}hler. However, Corollary \ref{cor3} yields, under growth assumption relative to Ricci flat K\"{a}hler metrics, the equivalence between Chern flatness and K\"{a}hlerness in noncompact context. 
\section{Local difference estimates and consequences}
	Let $(E,\bar\partial_E)$ be a holomorphic vector bundle over a K\"ahler manifold $(M,g)$. For a Hermitian metric $H$ on $E$, let $D_H=\partial_H+\bar\partial_E$ and $F_H=D_H^2$
	denote its Chern connection and Chern curvature respectively.
	Let $\nabla_H$ denote
	the connection induced by $D_H$ and the Chern connection
	$\nabla_g$. Given another Hermitian metric $K$, we define $h=K^{-1}H$. Let $u,v\in A^{0}(M,E)$, then
	\begin{equation}
		\begin{aligned}
			H(\pa_{H}u,v)
			&=\pa H(u,v)-H(u,\ol\pa_Ev)\\
			&=\pa K(h(u),v)-K(h(u),\ol\pa_Ev)\\
			&=K((\pa_K\circ h)(u),v)\\
			&=H((h^{-1}\circ\pa_K\circ h)(u),v)\\
			&=H((\pa_K+h^{-1}\pa_K h)(u),v).
		\end{aligned}
	\end{equation}
	Consequently, we have 
	\begin{equation}\label{connection-difference}
		D_H-D_K
		=h^{-1}\partial_Kh
		\in A^{1,0}(\End E).
	\end{equation}
	Since the Chern curvature is of type $(1,1)$, it follows
	\begin{equation}\label{curvature-difference}
		\begin{aligned}
			F_H&=F_K+\bar\partial_E(h^{-1}\partial_Kh).
		\end{aligned}
	\end{equation}
	\begin{proof}[\textup{\textbf{Proof of Theorem \ref{thm:main11}}}]
		Let $T_E=h^{-1}\pa_Kh$ and we take a holomorphic normal coordinate system centered at $x$.
		Below all computations are evaluated at $x$. 
	It holds
		\begin{equation}\label{proof311}
			\sqrt{-1}\Lambda\partial\bar\partial\big|T_E\big|^2
			=|\nabla_HT_E|^2
			+g^{\bar{j}i}
			\langle\nabla_{H,\frac{\partial}{\partial z^i}}
			\nabla_{H,\frac{\partial}{\partial\ol{z}^j}}T_E,T_E\rangle
			+g^{\bar{i}j}
			\langle T_E,\nabla_{H,\frac{\partial}{\partial\ol{z}^j}}
			\nabla_{H,\frac{\partial}{\partial z^i}}T_E\rangle.
		\end{equation}
		On the one hand, we deduce
		\begin{equation}
			\begin{aligned}\label{proof312}
				&\nabla_{H,\frac{\partial}{\partial z^i}}\nabla_{H,\frac{\partial}{\partial\ol{z}^j}}T_E(\frac{\partial}{\partial z^k})
				\\&=D_{H,\frac{\partial}{\partial z^i}}\big(\nabla_{H,\frac{\partial}{\partial\ol{z}^j}}T_E(\frac{\partial}{\partial z^k})\big)
				\\&=D_{H,\frac{\partial}{\partial z^i}}\big(F_K(\frac{\partial}{\partial z^k},\frac{\partial}{\partial\ol{z}^j})
				-F_H(\frac{\partial}{\partial z^k},\frac{\partial}{\partial\ol{z}^j})\big)
				\\&=D_{H,\frac{\partial}{\partial z^i}}\big(F_K(\frac{\partial}{\partial z^k},\frac{\partial}{\partial\ol{z}^j})\big)
				-D_{H,\frac{\partial}{\partial z^k}}\big(F_H(\frac{\partial}{\partial z^i},\frac{\partial}{\partial\ol{z}^j})\big)
				\\&=D_{K,\frac{\partial}{\partial z^i}}\big(F_K(\frac{\partial}{\partial z^k},\frac{\partial}{\partial\ol{z}^j})\big)
				-D_{H,\frac{\partial}{\partial z^k}}\big(F_H(\frac{\partial}{\partial z^i},\frac{\partial}{\partial\ol{z}^j})\big)+[T_E(\frac{\partial}{\partial z^i}),F_K(\frac{\partial}{\partial z^k},\frac{\partial}{\partial\ol{z}^j})]
				\\&=D_{K,\frac{\partial}{\partial z^k}}\big(F_K(\frac{\partial}{\partial z^i},\frac{\partial}{\partial\ol{z}^j})\big)
				-D_{H,\frac{\partial}{\partial z^k}}\big(F_H(\frac{\partial}{\partial z^i},\frac{\partial}{\partial\ol{z}^j})\big)+[T_E(\frac{\partial}{\partial z^i}),F_K(\frac{\partial}{\partial z^k},\frac{\partial}{\partial\ol{z}^j})],
			\end{aligned}
		\end{equation}
		where we have used $(\ref{curvature-difference})$ and the Bianchi identity.
		On the other hand, we have
		\begin{equation}
			\begin{aligned}\label{proof313}
				&\nabla_{H,\frac{\partial}{\partial\ol{z}^j}}\nabla_{H,\frac{\partial}{\partial z^i}}T_E(\frac{\partial}{\partial z^k})
				\\&=\nabla_{H,\frac{\partial}{\partial z^i}}\nabla_{H,\frac{\partial}{\partial\ol{z}^j}}T_E(\frac{\partial}{\partial z^k})
				-[F_H(\frac{\partial}{\partial z^i},\frac{\partial}{\partial\ol{z}^j}),T_E(\frac{\partial}{\partial z^k})]+T_E\big(F_{\nabla}(\frac{\partial}{\partial z^i},\frac{\partial}{\partial\ol{z}^j})\frac{\partial}{\partial z^k}\big)
				\\&=D_{K,\frac{\partial}{\partial z^k}}\big(F_K(\frac{\partial}{\partial z^i},\frac{\partial}{\partial\ol{z}^j})\big)
				-D_{H,\frac{\partial}{\partial z^k}}\big(F_H(\frac{\partial}{\partial z^i},\frac{\partial}{\partial\ol{z}^j})\big)+[T_E(\frac{\partial}{\partial z^i}),F_K(\frac{\partial}{\partial z^k},\frac{\partial}{\partial\ol{z}^j})]
				\\&-[F_H(\frac{\partial}{\partial z^i},\frac{\partial}{\partial\ol{z}^j}),T_E(\frac{\partial}{\partial z^k})]+T_E(F_{\nabla}\big(\frac{\partial}{\partial z^i},\frac{\partial}{\partial\ol{z}^j})\frac{\partial}{\partial z^k}\big).
			\end{aligned}
		\end{equation}
		It follows 
		\begin{equation}\label{proof314}
			\begin{aligned}
				&\sqrt{-1}\Lambda\partial\bar\partial|T_E|^2
				\\&=|\nabla_HT_E|^2+\langle D_{K}(\sqrt{-1}\Lambda F_K)
				-D_{H}(\sqrt{-1}\Lambda F_H),T_E\rangle
				\\&+g^{\bar{j}i}\langle [T_E(\frac{\partial}{\partial z^i}),F_K(\bullet,\frac{\partial}{\partial\ol{z}^j})],T_E\rangle
				+\langle T_E,D_{K}(\sqrt{-1}\Lambda F_K)-D_{H}(\sqrt{-1}\Lambda F_H)\rangle
				\\&+g^{\bar{i}j}\langle T_E, [T_E(\frac{\partial}{\partial z^i}),F_K(\bullet,\frac{\partial}{\partial\ol{z}^j})]\rangle
				-\langle T_E,[\sqrt{-1}\Lambda F_H,T_E]-T_E\circ\Ric^{(2)}\rangle
				\\&=|\nabla_HT_E|^2+2\operatorname{Re}\langle
				D_{K}(\sqrt{-1}\Lambda F_K)-D_{H}(\sqrt{-1}\Lambda F_H),T_E\rangle
				\\&+2\operatorname{Re}\langle g^{\bar{j}i} [T_E(\frac{\partial}{\partial z^i}),F_K(\bullet,\frac{\partial}{\partial\ol{z}^j})],T_E\rangle-\langle T_E,[\sqrt{-1}\Lambda F_H,T_E]-T_E\circ\Ric^{(2)}\rangle,
			\end{aligned}
		\end{equation}
	where $\Ric^{(2)}=\sqrt{-1}\Lambda F_\nabla$ is the second Chern-Ricci transformation.
		\par 
		Since $H$ is Hermitian-Einstein and $\det H=\det K$, we have 
	$[\sqrt{-1}\Lambda F_H,T_E]=0$ and
		\begin{equation}
			\begin{aligned}\label{proof315}
				&D_K(\sqrt{-1}\Lambda F_K)-
				D_H(\sqrt{-1}\Lambda F_H)\\
				&=D_K(\sqrt{-1}\Lambda F_K^\perp
				+\frac{\sqrt{-1}\tr\Lambda F_H}{\rank(E)}\Id_E)
				-D_H(\sqrt{-1}\Lambda F_H^\perp
				+\frac{\sqrt{-1}\tr\Lambda F_H}{\rank(E)}\Id_E)
				\\&=D_K(\sqrt{-1}\Lambda F_K^\perp).
			\end{aligned}
		\end{equation}
Since $\Ric^{(2)}$ coincides with the Ricci transformation, \eqref{proof314} yields 
		\begin{equation}\label{proof316}
			\begin{aligned}
				\sqrt{-1}\Lambda\pa\ol\pa|T_E|^2
				&\geq|\nabla_H T_E|^2-2|D_K(\Lambda F_K^\perp)||T_E|-4|F_K^{\perp}||T_E|^2-\kappa |T_E|^2.
			\end{aligned}
		\end{equation}
		Since $\sqrt{-1}\Lambda\pa\ol\pa|T_E|^2=2|T_E|\sqrt{-1}\Lambda\pa\ol\pa|T_E|+|\nabla|T_E||^2$, Kato inequality gives
		\begin{equation}\label{proof317}
			\begin{aligned}
				|T_E|\sqrt{-1}\Lambda\pa\ol\pa|T_E|
				\geq-\frac{\kappa}{2}|T_E|^2-2|F_K^{\perp}|
				|T_E|^2-|D_K(\Lambda F_K^\perp)||T_E|.
			\end{aligned}
		\end{equation}
		Thus at the point where $|T_E|\neq0$, we have
		\begin{equation}\label{proof318}
			\begin{aligned}
				\sqrt{-1}\Lambda\pa\ol\pa|T_E|
				&\geq-(\frac{\kappa}{2}+2|F_K^{\perp}|)|T_E|
				-|D_K(\Lambda F_K^\perp)|.
			\end{aligned}
		\end{equation}
		\par
		Let $\psi=\frac{|T_E|}{2A_R-\tr h^{-1}}$ and it holds
		\begin{equation}\begin{split}\label{proof319}
			\sqrt{-1}\Lambda\partial\bar\partial\psi
			&=\frac{\sqrt{-1}\Lambda\partial\bar\partial|T_E|}
			{2A_R-\tr h^{-1}}
			+|T_E|\sqrt{-1}\Lambda\big(\frac{\pa\ol\pa\tr h^{-1}}{(2A_R-\tr h^{-1})^2}+\frac{2\pa\tr h^{-1}\wedge\ol\pa\tr h^{-1}}{(2A_R-\tr h^{-1})^3}\big)
			\\&+\frac{2\operatorname{Re}\langle(2A_R-\tr h^{-1})\pa\psi-\psi\pa\tr h^{-1},\pa\tr h^{-1}\rangle}{(2A_R-\tr h^{-1})^2}
			\\&=
			\frac{\sqrt{-1}\Lambda\partial\bar\partial|T_E|}
			{2A_R-\tr h^{-1}}
			+\frac{\sqrt{-1}\Lambda\partial\bar\partial\tr h^{-1}}{2A_R-\tr h^{-1}}
			\psi
			+\frac{\langle\nabla\psi,\nabla\tr h^{-1}\rangle}
			{2A_R-\tr h^{-1}}.
		\end{split}\end{equation}
Applying $(\ref{curvature-difference})$ yields
\begin{equation}\begin{split}\label{proof3110}
\sqrt{-1}\Lambda\pa\ol\pa\tr h^{-1}
&=-\sqrt{-1}\Lambda\ol\pa\pa\langle h^{-1},\Id_E\rangle_K
\\&=\sqrt{-1}\Lambda\ol\pa\langle T_E,h^{-1}\rangle_K
\\&=\langle\sqrt{-1}\Lambda\ol\pa_ET_E,h^{-1}\rangle_K+\sqrt{-1}\Lambda\langle T_E,T_Eh^{-1}\rangle_K
\\&=-\langle\sqrt{-1}\Lambda F_K^\perp,h^{-1}\rangle_K+|T_Eh^{-\frac{1}{2}}|^2_{K,g},
	\end{split}\end{equation}
where we have used the Hermitian-Einstein condition and $\det H=\det K$. Thus $(\ref{proof318})$, $(\ref{proof319})$ and $(\ref{proof3110})$ show at the point where $|T_E|\neq0$, it holds
		\begin{equation}
			\begin{aligned}\label{proof3111}
				\sqrt{-1}\Lambda\pa\ol\pa\psi
				&\geq-(\frac{\kappa}{2}+2|F_K^\perp|
				+\frac{\langle\sqrt{-1}\Lambda F_K^\perp,h^{-1}\rangle_K}{2A_R-\tr h^{-1}}
				-\frac{|T_Eh^{-\frac{1}{2}}|_{K,g}^2}{2A_R-\tr h^{-1}})\psi
				\\
				&-\frac{|D_K(\Lambda F_K^\perp)|}
				{2A_R-\tr h^{-1}}+\frac{\langle\nabla\psi,\nabla\tr h^{-1}\rangle}
				{2A_R-\tr h^{-1}}.
			\end{aligned}
		\end{equation}
		Notice that
		\begin{equation}
			\begin{aligned}\label{proof3112}
				-\frac{\langle\sqrt{-1}\Lambda F_K^\perp, h^{-1}\rangle_K}{2A_R-\tr h^{-1}} 
				&\geq\frac{-\tr h^{-1}}{2A_R - \tr h^{-1}}|\Lambda F_K^\perp|_K
				\geq-|\Lambda F_K^\perp|_K.
			\end{aligned}
		\end{equation}
		Hence at the point where $|T_E|\neq0$, $(\ref{proof3111})$ gives
		\begin{equation}\label{proof3113}
			\begin{aligned}
				\sqrt{-1}\Lambda\pa\ol\pa\psi
				\geq&-(\frac{\kappa}{2}+2|F_K^\perp|+|\Lambda F_K^\perp|_K
			-\frac{|T_Eh^{-\frac{1}{2}}|_{K,g}^2}{2A_R-\tr h^{-1}})\psi\\
				&-\frac{|D_K(\Lambda F_K^\perp)|}
				{2A_R-\tr h^{-1}}+\frac{\langle\nabla\psi,\nabla\tr h^{-1}\rangle}
				{2A_R-\tr h^{-1}}.
			\end{aligned}
		\end{equation}
		\par We now introduce the auxiliary function on $\ol{B_R}$:
		\begin{equation}\label{proof3114}
			F=(R^2-\rho^2)\psi
			=\frac{|T_E|}{2A_R-\tr h^{-1}}
			(R^2-\rho^2),\ \rho=\dist(x,\bullet).
		\end{equation}
We may assume $T_E$ is not identically zero, then $F$ attains a positive maximum at some point $x_0\in B_R$. 
		If $x_0$ lies in the cut locus of $x$, Calabi's support function trick \cite{Ca1958} can be used to justify the argument below. All subsequent computations are considered at $x_0$. We have
		\begin{equation}\label{proof3115}
			\frac{\nabla\psi}{\psi}	=\frac{\nabla\rho^2}{R^2-\rho^2},\ 
			\frac{\sqrt{-1}\Lambda\pa\ol\pa\psi}{\psi}
			-\frac{\sqrt{-1}\Lambda\pa\ol\pa\rho^2}{R^2-\rho^2}
			-\frac{\langle\nabla\psi,\nabla\rho^2\rangle}{\psi(R^2-\rho^2)}\leq0,
		\end{equation}
and thus
		\begin{equation}\label{proof3116}
			0\geq\dfrac{\sqrt{-1}\Lambda\pa\ol\pa\psi}{\psi}
			-\dfrac{\sqrt{-1}\Lambda\pa\ol\pa\rho^2}{R^2-\rho^2}
			-\frac{|\nabla\rho^2|^2}{(R^2-\rho^2)^2}.
		\end{equation}
		It follows from \eqref{proof3113}, \eqref{proof3115} and \eqref{proof3116} that
		\begin{equation}\label{proof3117}
			\begin{aligned}
				0&\geq
				-(\frac{\kappa}{2}+2|F_K^\perp|
				+|\Lambda F_K^\perp|_K)
				+\frac{|T_Eh^{-\frac{1}{2}}|_{K,g}^2}{2A_R-\tr h^{-1}}
				\\
				&-\frac{|D_K(\Lambda F_K^\perp)|}
				{(2A_R-\tr h^{-1})\psi}
				+\frac{\langle\nabla\rho^2,\nabla\tr h^{-1}\rangle}
				{(2A_R-\tr h^{-1})(R^2-\rho^2)}\\
				&-\frac{\sqrt{-1}\Lambda\pa\ol\pa\rho^2}
				{R^2-\rho^2}-\frac{|\nabla\rho^2|^2}
				{(R^2-\rho^2)^2}.
			\end{aligned}
		\end{equation}
By the Laplacian comparison $\sqrt{-1}\Lambda\pa\ol\pa\rho^2
\leq C(1+\sqrt{\kappa}R)$ and $|\nabla\rho^2|^2\leq4R^2$, it gives 
		\begin{equation}\label{proof3118}
			\begin{aligned}
				\frac{|T_Eh^{-\frac{1}{2}}|_{K,g}^2}{2A_R-\tr h^{-1}}
				&\leq\frac{\kappa}{2}
				+2|F_K^\perp|+|\Lambda F_K^\perp|_K+\frac{|D_K(\Lambda F_K^\perp)|}{|T_E|}
				\\
				&+\frac{C(1+\sqrt{\kappa}R)}{R^2-\rho^2}
				+\frac{4R^2}{(R^2-\rho^2)^2}-\frac{\langle\nabla\rho^2,\nabla\tr h^{-1}\rangle}
				{(2A_R-\tr h^{-1})(R^2-\rho^2)}.
			\end{aligned}
		\end{equation}
For the last term on the right hand side, we note
\begin{equation}\label{proof3119A}
	\begin{aligned}
|\langle\nabla\rho^2,\nabla\tr h^{-1}\rangle|
&\leq\sqrt{2}|\nabla\rho^2||\pa\tr h^{-1}|
\\&=\sqrt{2}|\nabla\rho^2||\tr(T_E\circ h^{-1})|
\\&=\sqrt{2}|\nabla\rho^2||\langle T_E\circ h^{-1/2},h^{-1/2}\rangle_K|
\\&\leq\sqrt{2}|\nabla\rho^2||h^{-\frac{1}{2}}|_K||T_E\circ h^{-1/2}|_{K,g}
\\&=\sqrt{2\tr h^{-1}}|\nabla\rho^2||T_E\circ h^{-1/2}|_{K,g},
	\end{aligned}
\end{equation}
and then
		\begin{equation}\label{proof3119B}
			\begin{aligned}
-\frac{\langle\nabla\rho^2,\nabla\tr h^{-1}\rangle}
{(2A_R-\tr h^{-1})(R^2-\rho^2)}
				&\leq\frac{\sqrt{2\tr h^{-1}}|T_E\circ h^{-1/2}|_{K,g}|\nabla\rho^2|}
				{(2A_R-\tr h^{-1})(R^2-\rho^2)}\\
				&\leq\frac{|T_E\circ h^{-1/2}|_{K,g}^2}
				{2(2A_R-\tr h^{-1})}
				+\frac{\tr h^{-1}|\nabla\rho^2|^2}{(2A_R-\tr h^{-1})(R^2-\rho^2)^2}
				\\&\leq\frac{|T_E\circ h^{-1/2}|_{K,g}^2}
				{2(2A_R-\tr h^{-1})}
				+\frac{4R^2}{(R^2-\rho^2)^2}.
			\end{aligned}
		\end{equation}
Therefore \eqref{proof3118} and \eqref{proof3119B} imply
		\begin{equation}\label{proof3120}
			\begin{aligned}
		\frac{|T_E\circ h^{-1/2}|_{K,g}^2}
				{2A_R-\tr h^{-1}}
				&\leq\kappa
				+4|F_K^\perp|+2|\Lambda F_K^\perp|_K+\frac{2|D_K(\Lambda F_K^\perp)|}{|T_E|}
				\\
				&+\frac{C(1+\sqrt\kappa R)}{R^2-\rho^2}
				+\frac{16R^2}{(R^2-\rho^2)^2}.
			\end{aligned}
		\end{equation}
	For the term on the left hand side, we note 
			\begin{equation}\begin{split}\label{proof3121A}
|T_E|^2=|h^{\frac{1}{2}}\circ T_E\circ h^{-\frac{1}{2}}|_{K,g}^2
\leq\lambda_{R,\max}|T_E\circ h^{-\frac{1}{2}}|^2_{K,g},
	\end{split}\end{equation}
and then
		\begin{equation}\begin{split}\label{proof3121B}
\frac{|T_E\circ h^{-1/2}|_{K,g}^2}{2A_R-\tr h^{-1}}
&\geq\frac{|T_E|^2}{(2A_R-\tr h^{-1})\lambda_{R,\max}}		
\geq\frac{|T_E|^2}{2A_R\lambda_{R,\max}}.
	\end{split}\end{equation}
Multiplying \eqref{proof3120} by $(R^2-\rho^2)^2$ and using \eqref{proof3121B}, we have
	\begin{equation}\label{proof3122}
	\begin{aligned}
		(R^2-\rho^2)^2\frac{|T_E|^2}{2A_R\lambda_{R,\max}}
		&\leq\kappa R^4+C(1+\sqrt\kappa R)R^2+16R^2
		\\&+R^4(4|F_K^\perp|+2|\Lambda F_K^\perp|_K)+2R^4\frac{|D_K(\Lambda F_K^\perp)|}{|T_E|}
		\\&\leq CR^2(1+\sqrt\kappa R)^2
		\\&+R^4(4|F_K^\perp|+2|\Lambda F_K^\perp|_K)+2R^4\frac{|D_K(\Lambda F_K^\perp)|}{|T_E|}.
	\end{aligned}
\end{equation}
Taking the supremum of the curvature terms over $B_R$, we obtain
		\begin{equation}\label{proof3123}
			\begin{aligned}
				(R^2-\rho^2)^3|T_E|^3
				&\leq CA_R\lambda_{R,\max} R^2
				(1+\sqrt{\kappa}R)^2(R^2-\rho^2)|T_E|\\
				&+CA_R\lambda_{R,\max} R^4
				(\sup_{B_R}|F_K^\perp|
				+\sup_{B_R}|\Lambda F_K^\perp|_K)
				(R^2-\rho^2)|T_E|\\
				&+CA_R\lambda_{R,\max} R^6
				\sup\limits_{B_R}|
				D_K(\Lambda F_K^\perp)|.
			\end{aligned}
		\end{equation}
		By Young inequality, it holds		
		\begin{equation}
			\begin{aligned}\label{proof3124}
				(R^2-\rho^2)^3|T_E|^3
				&\leq C(A_R\lambda_{R,\max})^{3/2}R^3(1+\sqrt{\kappa}R)^3\\
				&+C(A_R\lambda_{R,\max})^{3/2}R^6
				(\sup\limits_{B_R}|F_K^\perp|+\sup\limits_{B_R}|\Lambda F_K^\perp|_K)^{3/2}\\
				&+CA_R\lambda_{R,\max}R^6
				\sup\limits_{B_R}|D_K(\Lambda F_K^\perp)|.
			\end{aligned}
		\end{equation}
Taking its cubic root yields
		\begin{equation}
			\begin{aligned}\label{proof3125}
				(R^2-\rho^2(x_0))|T_E(x_0)|
				&\leq C(A_R\lambda_{R,\max})^{1/2}
				R(1+\sqrt{\kappa}R)
				\\&+C(A_R\lambda_{R,\max})^{1/2}R^2(\sup_{B_R}|F_K^\perp|
+\sup_{B_R}|\Lambda F_K^\perp|_K)^{1/2}
\\&+C(A_R\lambda_{R,\max})^{1/3}R^2
\sup_{B_R}|D_K(\Lambda F_K^\perp)|^{1/3}.
			\end{aligned}
		\end{equation}	
Finally the conclusion follows from \eqref{proof3125} and
	\begin{equation}
	\begin{aligned}\label{proof3126}
		\sup_{B_{R/2}}|T_E|
		&=\sup_{B_{R/2}}\frac{2A_R-\tr h^{-1}}{R^2-\rho^2}F
		\\&\leq\frac{CA_R}{R^2}\frac{R^2-\rho^2(x_0)}{2A_R-\tr h^{-1}(x_0)}|T_E|(x_0)
		\\&\leq\frac{C}{R^2}(R^2-\rho^2(x_0))|T_E(x_0)|.
	\end{aligned}
\end{equation}
	\end{proof}
	\begin{proof}[\textup{\textbf{Proof of Theorem \ref{thm:main22}}}]
		We use notations in the proof of Theorem \ref{thm:main11}. We note
		\begin{equation}\label{proof411}
			\begin{aligned}
				\sqrt{-1}\Lambda\pa_K[F_K,h](\frac{\pa}{\pa z^k})
				&=g^{\ol{j}i}(\iota_{\frac{\pa}{\pa\ol{z}^j}}\iota_{\frac{\pa}{\pa z^i}}\pa_K[F_K,h])(\frac{\pa}{\pa z^k})
				\\&=g^{\ol{j}i}\pa_K[F_K,h](\frac{\pa}{\pa z^i},\frac{\pa}{\pa\ol{z}^j},\frac{\pa}{\pa z^k})
				\\&=-g^{\ol{j}i}\pa_K[F_K,h](\frac{\pa}{\pa z^i},\frac{\pa}{\pa z^k},\frac{\pa}{\pa\ol{z}^j}).
			\end{aligned}
		\end{equation}
		 By Bianchi identity, we have
		\begin{equation}\label{proof412}
			\begin{aligned}
				&\pa_K[F_K,h](\frac{\pa}{\pa z^i},\frac{\pa}{\pa z^k},\frac{\pa}{\pa\ol{z}^j})
				\\&=(\nabla_{K,\frac{\partial}{\partial z^i}}[F_K,h])
				(\frac{\partial}{\partial z^k},\frac{\partial}{\partial\bar z^j})
				-(\nabla_{K,\frac{\partial}{\partial z^k}}[F_K,h])
				(\frac{\partial}{\partial z^i},\frac{\partial}{\partial\bar z^j})\\
				&=[\nabla_{K,\frac{\partial}{\partial z^i}}F_K
				(\frac{\partial}{\partial z^k},\frac{\partial}{\partial\bar z^j})
				-\nabla_{K,\frac{\partial}{\partial z^k}}F_K
				(\frac{\partial}{\partial z^i},\frac{\partial}{\partial\bar z^j}),h]\\
				&+[F_K(\frac{\partial}{\partial z^k},\frac{\partial}{\partial\bar z^j}),
				D_{K,\frac{\partial}{\partial z^i}}h]
				-[F_K(\frac{\partial}{\partial z^i},\frac{\partial}{\partial\bar z^j}),
				D_{K,\frac{\partial}{\partial z^k}}h]
				\\
				&=[F_K(\frac{\partial}{\partial z^k},\frac{\partial}{\partial\bar z^j}),
				hT_E(\frac{\partial}{\partial z^i})]
				-[F_K\left(\frac{\partial}{\partial z^i},\frac{\partial}{\partial\bar z^j}\right),
				hT_E\left(\frac{\partial}{\partial z^k}\right)]\\
				&=[F_K(\frac{\partial}{\partial z^k},\frac{\partial}{\partial\bar z^j}),h]
				T_E(\frac{\partial}{\partial z^i})
				+h[F_K(\frac{\partial}{\partial z^k},\frac{\partial}{\partial\bar z^j}),
				T_E(\frac{\partial}{\partial z^i})]\\
				&-[F_K(\frac{\partial}{\partial z^i},\frac{\partial}{\partial\bar z^j}),h]
				T_E(\frac{\partial}{\partial z^k})
				-h[F_K(\frac{\partial}{\partial z^i},\frac{\partial}{\partial\bar z^j}),
				T_E(\frac{\partial}{\partial z^k})].
			\end{aligned}
		\end{equation}
Hence 
	\begin{equation}\label{proof413}
	\begin{aligned}
-\sqrt{-1}\Lambda\pa_K[F_K,h]
&=g^{\ol{j}i}[F_K(\bullet,\frac{\partial}{\partial\bar z^j}),h]
T_E(\frac{\partial}{\partial z^i})
+g^{\ol{j}i}h[F_K(\bullet,\frac{\partial}{\partial\bar z^j}),
T_E(\frac{\partial}{\partial z^i})]\\
&-[\sqrt{-1}\Lambda F_K,h]
T_E
-h[\sqrt{-1}\Lambda F_K,
T_E],
	\end{aligned}
\end{equation}
and equivalently,
		\begin{equation}\label{proof414}
			\begin{aligned}
g^{\ol{j}i}[F_K(\bullet,\frac{\partial}{\partial\ol{z}^j}),
T_E(\frac{\partial}{\partial z^i})]
&=-h^{-1}\sqrt{-1}\Lambda\pa_K[F_K,h]
-g^{\ol{j}i}h^{-1}[F_K(\bullet,\frac{\partial}{\partial\bar z^j}),h]
T_E(\frac{\partial}{\partial z^i})
\\&+h^{-1}[\sqrt{-1}\Lambda F_K,h]T_E+[\sqrt{-1}\Lambda F_K,T_E].
			\end{aligned}
		\end{equation}
It is easy to see $\Lambda F_K^{\ast H}=-h^{-1}\circ\Lambda F_K\circ h$ and thus	
\begin{equation}\begin{split}\label{proof415}
|h^{-1}[\Lambda F_K,h]|^2
&=|-\Lambda F_K^{\ast H}-\Lambda F_K|^2
\\&=2|\Lambda F_K^\perp|^2-2|\Lambda F_K^\perp|^2_K
\\&\leq2|\Lambda F_K^\perp|^2.
\end{split}\end{equation}
		It follows from \eqref{proof414} and \eqref{proof415}  that
		\begin{equation}\label{proof416}
			\begin{aligned}
				2\operatorname{Re}g^{i\bar j}
				\langle[T_E(\partial_i),F_K^\perp(\bullet,\bar\partial_j)],
				T_E\rangle
				&\geq-C|h^{-1}
				\Lambda\partial_K[F_K,h]||T_E|\\
				&-C|h^{-1}[F_K,h]||T_E|^2
				-C|\Lambda F_K^\perp||T_E|^2.
			\end{aligned}
		\end{equation}	
Since $H$ is Hermitian-Einstein and $\det H=\det K$, we have	
\begin{equation}\label{proof417}
	\begin{aligned}
\langle D_H(\sqrt{-1}\Lambda F_H),T_E\rangle=0,
	\end{aligned}
\end{equation}	\begin{equation}\label{proof418}
\begin{aligned}
&\operatorname{Re}\langle D_K(\sqrt{-1}\Lambda F_K),T_E\rangle
\geq-|D_K(\Lambda F_K^\perp)||T_E|.
\end{aligned}
\end{equation}
		Combining \eqref{proof314}, \eqref{proof416}, \eqref{proof417} and \eqref{proof418}, we obtain
		\begin{equation}\label{proof419}
			\begin{aligned}
				\sqrt{-1}\Lambda\pa\ol\pa|T_E|^2
				&\geq|\nabla_HT_E|^2-\kappa|T_E|^2-2|D_K(\Lambda F_K^\perp)||T_E|
				\\&-C|h^{-1}
				\Lambda D_K[F_K,h]||T_E|
				-C|h^{-1}[F_K,h]||T_E|^2
				-C|\Lambda F_K^\perp||T_E|^2.
			\end{aligned}
		\end{equation}
Applying Kato inequality then implies at the point where $|T_E|>0$, it holds
		\begin{equation}\label{proof4110}
			\begin{aligned}
				\sqrt{-1}\Lambda\pa\ol\pa|T_E|
				\geq&-(\frac{\kappa}{2}+C|h^{-1}[F_K,h]|+C|\Lambda F_K^\perp|)|T_E|
				\\&-|D_K(\Lambda F_K^\perp)|
				-C|h^{-1}
				\Lambda D_K[F_K,h]|.
			\end{aligned}
		\end{equation}
Running the same argument 
		as in the proof of Theorem~\ref{thm:main11}, we finish the proof.
	\end{proof}

	\section{Applications on holomorphic tangent bundles}
	Let $M$ be a complex manifold with complex dimension $n$ and complex structure $J$. In the presence of $J$, the tangent bundle $TM$ is a complex vector bundle. Since $J^2=-\Id$, we have
	\begin{equation}\begin{split}
			TM\otimes\mathbb{C}=T^{1,0}M\oplus T^{0,1}M,
	\end{split}\end{equation}
	where $T^{1,0}M$ and $T^{0,1}M$ are the eigenspaces of eigenvalues $\sqrt{-1}$ and $-\sqrt{-1}$ respectively. In fact, for any $p\in M$ we can find vectors $\{e_i,Je_i\}_{i=1}^{n}$ forming a basis of $T_pM$ and then $T^{1,0}_pM$ (resp. $T^{0,1}_pM$) is generated by $\{e_i-\sqrt{-1}Je_i\}_{i=1}^n$ (resp. $\{e_i+\sqrt{-1}Je_i\}_{i=1}^n$). We often identify $TM$ with $T^{1,0}M$ via
	\begin{equation}\begin{split}\label{identification}
			e_i\leftrightarrow\frac{1}{\sqrt{2}}(e_i-\sqrt{-1}Je_i),\  Je_i\leftrightarrow\frac{\sqrt{-1}}{\sqrt{2}}(e_i-\sqrt{-1}Je_i),\
			i=1,...,n.
	\end{split}\end{equation}
It is easy to see $T^{1,0}M$ is a complex vector bundle whose transition functions are holomorphic. On the one hand, via the identification $(\ref{identification})$, the complex vector bundle $TM$ was equipped with a holomorphic structure $\ol\pa_{TM}$ coming from that of $T^{1,0}M$. On the other hand, given a Hermitian metric $g$ on $M$, we can equip $TM$ with a Hermitian metric
	\begin{equation}\begin{split}\label{Hg}
H_g=g-\sqrt{-1}\omega_g,\ \omega_g(\bullet,\bullet)=g(J(\bullet),\bullet).
	\end{split}\end{equation}
In summary, $TM$ possesses a holomorphic structure $\ol\pa_{TM}$ and $H_g$. It is known that $g$ is K\"{a}hler if and only if $D_{H_g}=\nabla_g$, where $D_{H_g}$ is the Chern connection on $TM$.
\begin{proof}[\textup{\textbf{Proof of Theorem \ref{thm:main3}}}]
By assumption, we have 
\begin{equation}\begin{split}
\log c_R\leq f\triangleq\frac{1}{n}\log\frac{\omega^n}{\omega_0^n}\leq\log C_R.
\end{split}\end{equation}
Since $\omega$ and $\omega_0$ are Ricci flat, we have $\pa\ol\pa f=0$ and thus $f$ is harmonic with respect to $g_0$. The inequality $c_R\omega_0|_{B_R}\leq\omega|_{B_R}$ gives $c_R\leq B$ for a constant $B$ independent of $R$. Consequently, 
\begin{equation}\begin{split}
\frac{C_R}{c_R}
	\le \frac{BC_R}{c_R^2}=o(R^2).
\end{split}\end{equation}
In particular, $0\leq\log\frac{C_R}{c_R}\leq2\log R$ for large $R$.
Then Theorem \ref{thm:CY} implies
\begin{equation}\begin{split}
	\sup_{B_{R/2}}|df|_{g_0}
	\leq \frac{C}{R}\log R
	=o(1).
\end{split}\end{equation}
It follows $f=a$ for some constant $a$ and $\omega^n=e^{an}\omega_0^n$. Let $g$ be the Hermitian metric corresponding to $\omega$, then it is easy to see 
\begin{equation}\begin{split}
		\det H_g=\det H_{e^ag_0}.
\end{split}\end{equation}
Let $B^g_R=\{z\in\mathbb{C}^n|\dist_g(0,z)<R\}$, it follows
\begin{equation}\begin{split}
	B_{\sqrt{c_R}R}^{g}\subset B_R,\ \sqrt{c_R}R\geq
	\sqrt{\frac{c_R^2R^2}{C_R}}\rightarrow\infty.
\end{split}\end{equation}
Hence $(\mathbb{C}^n,g)$ is complete and it holds
\begin{equation}\begin{split}
		c_R\omega_0|_{B_{\sqrt{c_R}R}^{g}}\leq\omega|_{B_{\sqrt{c_R}R}^{g}}\leq C_R\omega_0|_{B_{\sqrt{c_R}R}^{g}},\ \frac{C_R}{c_R}=o(c_RR^2).
\end{split}\end{equation}
Therefore Corollary \ref{cor1} applies to conclude $D_{H_g}=D_{H_{g_0}}$ and thus $\nabla_g=\nabla_{g_0}$. Since $\nabla_g$ is the Levi-Civita connection with respect to $g$, we have
\begin{equation}\begin{split}
\nabla_{g_0}\omega=\nabla_{g}\omega=0.	
\end{split}\end{equation}
\end{proof}
\begin{proof}[\textup{\textbf{Proof of Theorem \ref{thm4}}}]
	Let $\omega$ and $\tilde\omega$ be fundamental forms of $g$ and $\tilde{g}$ respectively. As in the proof of Theorem \ref{thm:main3}, we have $\frac{1}{n}\log\frac{\tilde\omega^n}{\omega^n}=a$ for some constant $a$ and thus 
	\begin{equation}\begin{split}
		\det H_{e^ag}=\det H_{\tilde{g}}.
	\end{split}\end{equation}
	Therefore Corollary \ref{cor1} applies to conclude $D_{H_g}=D_{H_{\tilde{g}}}$ and thus $\nabla_g=\nabla_{\tilde{g}}$.
\end{proof}
\section*{Appendix: a counterexample}
	\par
	In this section, we give an explicit example showing that on complete noncompact K\"{a}hler manifolds with nonnegative Ricci curvature, two mutually bounded Hermitian-Einstein metrics with same determinant may have different Chern connections.
\par Let $E=\mathbb C^2\times\mathbb C^2$ be the trivial holomorphic vector bundle
		over $\mathbb C^2$ equipped with the standard K\"ahler metric. For a constant $a>0$, we consider 
		\begin{equation}\label{counterexample-Ha}
			H_a
			=
			\frac{a+|z|^2}{a}\Id_E
			-\frac{2a+|z|^2}{a(a+|z|^2)}\ol{z}z^T,
		\end{equation}
where $z=(z^1,z^2)^{T}$. For $z\neq0$ and $P=\frac{\ol{z}z^T}{|z|^2}$, we have
		\begin{equation}\label{counterexample-metricPro}
			H_a
			=
			\frac{a}{a+|z|^2}P
			+
			\frac{a+|z|^2}{a}(\Id_E-P).
		\end{equation}
As a Hermitian matrix function, it is easy to see it has eigenvalues 
\begin{equation}\begin{split}
\frac{a}{a+|z|^2},\ \frac{a+|z|^2}{a},
\end{split}\end{equation}	
and eigenvectors $\ol{z}$, $\ol{z}^\perp$ respectively, where $\ol{z}^\perp$ is orthogonal to $\ol{z}$. Moreover, for any positive constants $a\neq b$, we have
\begin{equation}\begin{split}
\min\{\frac{a}{b},\frac{b}{a}\}H_a\leq H_b\leq\max\{\frac{a}{b},\frac{b}{a}\}H_a.
	\end{split}
\end{equation} 
Therefore $H_a$, $H_b$ are two mutually bounded Hermitian metrics and $\det H_a=\det H_b$. Let $(e_1,e_2)$ be the standard global frame of $E$ and $\pa_{H_a}(e_1,e_2)=(e_1,e_2)A_{H_a}$, then we have
\begin{equation}\begin{split}
A_{H_a}=\ol{H}_a^{-1}\pa\ol{H}_a,\ F_{H_a}=\ol\pa A_{H_a}.	
\end{split}\end{equation}
A straightforward computation shows
\begin{equation}\label{counterexample-connection-form}
\begin{aligned}
A_{H_a}=\sum\limits_{i=1}^2\frac{\ol{z}^idz^i}{a+|z|^2}\Id_E-\frac{2a+|z|^2}{(a+|z|^2)^2}dz\ol{z}^T-\sum\limits_{i=1}^2\frac{\ol{z}^idz^i}{(a+|z|^2)^2}z\ol{z}^T,
\end{aligned}
		\end{equation}
		\begin{equation}\label{counterexample-connection-form}
		\begin{aligned}
F_{H_a}
&=\big(-\sum\limits_{i=1}^2\frac{dz^i\wedge d\ol{z}^i}{a+|z|^2}
+\sum\limits_{i,j=1}^2\frac{\ol{z}^iz^jdz^i\wedge d\ol{z}^j}{(a+|z|^2)^2}\big)\Id_E
\\&+\frac{2a+|z|^2}{(a+|z|^2)^2}dz\wedge d\ol{z}^T-\sum\limits_{i=1}^2\frac{(3a+|z|^2)z^i}{(a+|z|^2)^3}(dz\wedge d\ol{z}^i)\ol{z}^T
\\&+\sum\limits_{i=1}^2\frac{dz^i\wedge d\ol{z}^i}{(a+|z|^2)^2}z\ol{z}^T
-\sum\limits_{i,j=1}^2\frac{2\ol{z}^iz^jdz^i\wedge d\ol{z}^j}{(a+|z|^2)^3}z\ol{z}^T
+\sum\limits_{i=1}^2\frac{z\ol{z}^idz^i\wedge d\ol{z}^T}{(a+|z|^2)^2}.
		\end{aligned}
	\end{equation}
Therefore $\sqrt{-1}\Lambda F_{H_a}=0$, both $H_a,H_b$ are Hermitian-Einstein but $D_{H_a}\neq D_{H_b}$.
\addtocontents{toc}{\protect\setcounter{tocdepth}{-1}}
\section*{AI declaration}
\addtocontents{toc}{\protect\setcounter{tocdepth}{2}}
ChatGPT was solely used to identify a counterexample in the Appendix and correct typographical errors.
\addtocontents{toc}{\protect\setcounter{tocdepth}{-1}}
\section*{Acknowledgements}
\addtocontents{toc}{\protect\setcounter{tocdepth}{2}}
We thank Chao Li for useful discussions on the topic of this paper.


\begin{thebibliography}{1}
		\footnotesize
		\bibitem[ADM]{ADM1996}
		Vestislav Apostolov, Johann Davidov and Oleg Mu\v{s}karov,  Compact self-dual Hermitian surfaces. Transactions of the American Mathematical Society. 348 (1996), no. 8, 3051-3063.
		\bibitem[Bal]{Bal1985}
        Andrew Balas, Compact Hermitian manifolds of constant holomorphic sectional curvature. Mathematische Zeitschrift. 189 (1985), no. 2, 193-210.
        \bibitem[BG]{BG1985} 
        Andrew Balas and Paul Gauduchon, Any Hermitian metric of constant nonpositive (Hermitian) holomorphic sectional curvature on a compact complex surface is K\"{a}hler. Mathematische Zeitschrift. 190 (1985), no. 1, 39-43.
		\bibitem[Ban]{Ba1993}
		Shigetoshi Bando, Einstein-Hermitian metrics on noncompact K\"{a}hler manifolds. 
		Einstein metrics and Yang-Mills connections (Sanda, 1990), 27-33, 
		Lecture Notes in Pure and Appl. Math. 145, 
		Dekker, New York, 1993.
		\bibitem[Bo]{Bo1958}
        William M Boothby, Hermitian manifolds with zero curvature. Michigan Mathematical Journal. 5 (1958), 229-233.
	\bibitem[Ca]{Ca1958}	
		Eugenio Calabi, An extension of E. Hopf's maximum principle with an application to Riemannian geometry. Duke Mathematical Journal. 25 (1958), 45-56. 
		
		\bibitem[CZ]{CZ2026}
		Shuwen Chen and Fangyang Zheng, Hermitian threefolds with constant holomorphic sectional curvature. arXiv: 2609.39282, 2026.
		
		\bibitem[CY]{ChengYau75}
		Shiu Yuen Cheng and Shing-Tung Yau, Differential equations on Riemannian manifolds and their geometric applications. 
		Communications on Pure and Applied Mathematics. 
		28 (1975), no. 3, 333-354.
		
		
		
		\bibitem[Do1]{Donaldson85}
		Simon Kirwan Donaldson, Anti self-dual Yang-Mills connections over complex algebraic surfaces and stable vector bundles. 
		Proceedings of the London Mathematical Society. 
		50 (1985), no. 1, 1-26.
		
		\bibitem[Do2]{Donaldson87}
		Simon Kirwan Donaldson, Infinite determinants, stable bundles and curvature. 
		Duke Mathematical Journal. 54 (1987), no. 1, 231-247.
		
		
		
		%
		
		
		
		
		
		
		
		
		
		
		
		
		
		
		
		
		
		%
		%
		
		\bibitem[MR]{MR1984}
		Vikram Bhagvandas Mehta and Annamalai Ramanathan, 
		Restriction of stable sheaves and representations of the fundamental group. 
		Inventiones Mathematicae. 77 (1984), no. 1, 163-172.
		
		\bibitem[Mo]{MochizukiHE}
		Takuro Mochizuki, Kobayashi--Hitchin correspondence for analytically stable bundles. 
		Transactions of the American Mathematical Society. 
		373 (2020), no. 1, 551-596.
		
		\bibitem[MS]{MS2010}
		Ovidiu Munteanu and Natasa Sesum, The Poisson equation on complete manifolds with positive spectrum and applications. Advances in Mathematics. 223 (2010), no. 1, 198-219.
		
		\bibitem[NS]{NS1965}
		Mudumbai Seshachalu Narasimhan and Conjeeveram Srirangachari Seshadri, 
		Stable and unitary vector bundles on a compact Riemann surface. 
		Annals of Mathematics. 82 (1965), 540-567.
		
		\bibitem[Ni]{Ni2002}
		Lei Ni, The Poisson equation and Hermitian-Einstein metrics on holomorphic vector bundles over complete noncompact K\"{a}hler manifolds. 
		Indiana University Mathematics Journal. 
		51 (2002), no. 3, 679-704.
		
		\bibitem[NR]{NiRen}
		Lei Ni and Huaiyu Ren, Hermitian-Einstein metrics for vector bundles on complete K\"{a}hler manifolds. 
		Transactions of the American Mathematical Society. 
		353 (2001), no. 2, 441-456.
		
		
		
		
		
		\bibitem[QT]{QT2026}
		Xiangsen Qin and Yuanhong Tian, Chern Holomorphic Sectional Curvature: Canonical Ampleness, Rigidity, and Flatness Counterexamples. arXiv:2609.25762v7, 2026. 
		\bibitem[RS]{RS1984}
		 Dieter Riebesehl and Friedmar Schulz, A priori estimates and a Liouville theorem for complex Monge-Amp\`{e}re equations. Mathematische Zeitschrift. 186 (1984), no. 1, 57-66.
		
		
		\bibitem[Si]{Si1988}
		Carlos Tschudi Simpson, Constructing variations of Hodge structure using Yang-Mills theory and applications to uniformization. 
		Journal of the American Mathematical Society. 
		1 (1988), no. 4, 867-918.
		
		
		
		
		\bibitem[UY]{UhlenbeckYau86}
		Karen Keskulla Uhlenbeck and Shing-Tung Yau, On the existence of Hermitian-Yang-Mills connections in stable vector bundles. Frontiers of the mathematical sciences: 1985 (New York, 1985). Communications on Pure and Applied Mathematics. 
		39 (1986), no. S, suppl., S257-S293.
		
		
		
		\bibitem[WZ]{WuZhangHE}
		Di Wu and Xi Zhang, Hermitian-Einstein equations on noncompact manifolds, 
		Proceedings of the London Mathematical Society. 
		133 (2026), no. 3, Paper No. e70216, 37 pp.
		
		\bibitem[Ya]{Yau75}
		Shing-Tung Yau, Harmonic functions on complete Riemannian manifolds. 
		Communications on Pure and Applied Mathematics. 
		28 (1975), 201-228.
		
		\bibitem[Zh]{ZhangJ}
		Junsheng Zhang, Hermitian-Yang-Mills connections on some complete non-compact K\"{a}hler manifolds. Mathematische Annalen. 390 (2024), no. 3, 4535-4575.
		
	\end{thebibliography}
\end{document}